\documentclass[11pt]{article}      %{article} was 12pt latex e
\usepackage{amssymb}
\usepackage{amscd}

\newtheorem{theorem}{Theorem}[section]
\newtheorem{prop}[theorem]{Proposition}
\newtheorem{defn}[theorem]{Definition}
\newtheorem{lemma}[theorem]{Lemma}
\newtheorem{coro}[theorem]{Corollary}
\newtheorem{prop-def}{Proposition-Definition}[section]

\newtheorem{remark}{Remark}[section]
\newtheorem{example}{Example}[section]

\newcommand{\nc}{\newcommand}

\nc{\bin}{\binc}
\nc{\binc}[2]{\left (\!\! \begin{array}{c} \scs{#1}\\
    \scs{#2} \end{array}\!\!\right )}  %binomial coeff 
\nc{\bincc}{\binc}
\nc{\bs}{\bar{S}}
\nc{\la}{\longrightarrow}
\nc{\rar}{\rightarrow}
\nc{\dar}{\downarrow}
\nc{\dap}[1]{\downarrow \rlap{$\scriptstyle{#1}$}}
\nc{\defeq}{\stackrel{\rm def}{=}}
\nc{\dis}[1]{\displaystyle{#1}}
\nc{\dotcup}{\ \displaystyle{\bigcup^\bullet}\ }
\nc{\hcm}{\ \hat{,}\ }
\nc{\hts}{\hat{\otimes}}
\nc{\hcirc}{\hat{\circ}}
\nc{\lleft}{[}
\nc{\lright}{]}
\nc{\curlyl}{\left \{ \begin{array}{c} {} \\ {} \end{array}
    \right .  \!\!\!\!\!\!\!} 
\nc{\curlyr}{ \!\!\!\!\!\!\!
    \left . \begin{array}{c} {} \\ {} \end{array}
    \right \} }
\nc{\longmid}{\left | \begin{array}{c} {} \\ {} \end{array}
    \right . \!\!\!\!\!\!\!}    %long vertical line
\nc{\ora}[1]{\stackrel{#1}{\rar}}
\nc{\ola}[1]{\stackrel{#1}{\la}}%${\Bbb Z}$
\nc{\scs}[1]{\scriptstyle{#1}}
\nc{\mrm}[1]{{\rm #1}}
\nc{\margin}[1]{\marginpar{}}   %{\rm #1}}
\nc{\dirlim}{\displaystyle{\lim_{\longrightarrow}}\,}
\nc{\invlim}{\displaystyle{\lim_{\longleftarrow}}\,}
\nc{\dislim}[1]{\displaystyle{\lim_{#1}}}
\nc{\colim}{\mrm{colim}}
\nc{\mvp}{\vspace{0.3cm}}
\nc{\tk}{^{(k)}}
\nc{\tp}{^\prime}
\nc{\ttp}{^{\prime\prime}}
\nc{\svp}{\vspace{2cm}}
\nc{\vp}{\vspace{8cm}}
\nc{\proofend}{ $\blacksquare$ \vspace{0.3cm}}
\nc{\modg}[1]{\!<\!\!{#1}\!\!>}
\nc{\intg}[1]{F_C(#1)}
\nc{\lmodg}{\!<\!\!}
\nc{\rmodg}{\!\!>\!}
\nc{\cpi}{\widehat{\Pi}}
\nc{\sha}{{\mbox{\cyr X}}}  %used to be \cyr
\nc{\shpr}{\diamond}   %Shuffle product
\nc{\labs}{\mid\!}
\nc{\rabs}{\!\mid}
\nc{\hatsha}{\widehat{\sha}}

\nc{\ann}{\mrm{ann}}
\nc{\Aut}{\mrm{Aut}}
\nc{\can}{\mrm{can}}
\nc{\Cont}{\mrm{Cont}}
\nc{\rchar}{\mrm{char}}
\nc{\cok}{\mrm{coker}}
\nc{\dtf}{{R-{\rm tf}}}
\nc{\dtor}{{R-{\rm tor}}}

\nc{\Div}{{\mrm Div}}
\nc{\End}{\mrm{End}}
\nc{\Ext}{\mrm{Ext}}
\nc{\Fil}{\mrm{Fil}}
\nc{\Fr}{\mrm{Fr}}
\nc{\Frob}{\mrm{Frob}}
\nc{\Gal}{\mrm{Gal}}
\nc{\GL}{\mrm{GL}}
\nc{\Hom}{\mrm{Hom}}
\nc{\hsr}{\mrm{H}}
\nc{\hpol}{\mrm{HP}}
\nc{\id}{\mrm{id}}
\nc{\im}{\mrm{im}}
\nc{\incl}{\mrm{incl}}
\nc{\length}{\mrm{length}}
\nc{\mchar}{\rm char}
\nc{\mpart}{\mrm{part}}
\nc{\ql}{{\QQ_\ell}}
\nc{\qp}{{\QQ_p}}
\nc{\rank}{\mrm{rank}}
\nc{\rcot}{\mrm{cot}}
\nc{\rdef}{\mrm{def}}
\nc{\rdiv}{{\rm div}}
\nc{\rtf}{{\rm tf}}
\nc{\rtor}{{\rm tor}}
\nc{\res}{\mrm{res}}
\nc{\SL}{\mrm{SL}}
\nc{\Spec}{\mrm{Spec}}
\nc{\tor}{\mrm{tor}}
\nc{\Tr}{\mrm{Tr}}
\nc{\tr}{\mrm{tr}}

\nc{\ab}{\mathbf{Ab}}
\nc{\Alg}{\mathbf{Alg}}
\nc{\Bax}{\mathbf{Bax}}
\nc{\bfk}{{\bf k}}
\nc{\bfone}{{\bf 1}}
\nc{\detail}{\marginpar{\bf More detail}
    \noindent{\bf Need more detail!}
    \svp}
\nc{\Diff}{\mathbf{Diff}}   
\nc{\gap}[1]{\marginpar{\bf Gap}\noindent{\bf #1}
    \svp}
\nc{\FMod}{\mathbf{FMod}}
\nc{\Int}{\mathbf{Int}}
\nc{\Mon}{\mathbf{Mon}}
\nc{\proof}{\noindent{\bf Proof: }}
\nc{\remarks}{\noindent{\bf Remarks: }}
\nc{\Rep}{\mathbf{Rep}}
\nc{\Rings}{\mathbf{Rings}}
\nc{\Sets}{\mathbf{Sets}}
\nc{\bill}[1]{\marginpar{\bf To Bill}\noindent{\bf To Bill:}
    {\tt #1}\\ }
\nc{\li}[1]{\marginpar{\bf To Li}\noindent{\bf To Li:}
    {\tt #1}\\ }

\nc{\BA}{{\mathbb A}}
\nc{\CC}{{\mathbb C}}
\nc{\DD}{{\mathbb D}}
\nc{\EE}{{\mathbb E}}
\nc{\FF}{{\mathbb F}}
\nc{\GG}{{\mathbb G}}
\nc{\HH}{{\mathbb H}}
\nc{\LL}{{\mathbb L}}
\nc{\NN}{{\mathbb N}}
\nc{\QQ}{{\mathbb Q}}
\nc{\RR}{{\mathbb R}}
\nc{\TT}{{\mathbb T}}
\nc{\VV}{{\mathbb V}}
\nc{\ZZ}{{\mathbb Z}}

\nc{\cala}{{\cal A}}
\nc{\calc}{{\cal C}}
\nc{\cald}{\mathcal{D}}
\nc{\cale}{{\cal E}}
\nc{\calf}{{\cal F}}
\nc{\calg}{{\cal G}}
\nc{\calh}{{\cal H}}
\nc{\cali}{{\cal I}}
\nc{\call}{{\cal L}}
\nc{\calm}{{\cal M}}
\nc{\caln}{{\cal N}}
\nc{\calo}{{\cal O}}
\nc{\calp}{{\cal P}}
\nc{\calr}{{\cal R}}
\nc{\cals}{{\cal S}}
\nc{\calt}{{\cal T}}
\nc{\calw}{{\cal W}}
\nc{\calx}{{\cal X}}
\nc{\caly}{{\cal Y}}
\nc{\CA}{\mathcal{A}}

\nc{\fraka}{{\mathfrak a}}
\nc{\frakA}{{\mathfrak A}}
\nc{\frakB}{{\mathfrak B}}
\nc{\frakm}{{\mathfrak m}}
\nc{\frakp}{{\mathfrak p}}
\nc{\frakS}{{\mathfrak S}}

\nc{\bfn}{\NN}

\font\cyr=wncyr10

\title{Ascending Chain Conditions in Free Baxter Algebras
\thanks{
The author is supported in part by NSF grant
    \#DMS 97-96122. 
MSC Numbers: 13E05, 16W99, 47C05.
Keywords: Baxter algebra, ascending chain condition}
\\
}
\author{
Li Guo\\
Department of Mathematics and Computer Science\\
Rutgers University\\
Newark, NJ 07102 \\
USA
}

\date{}

\begin{document}
\maketitle

\begin{abstract}
In this paper we study ascending chain conditions 
in a free Baxter algebra by making use of explicit constructions 
of free Baxter algebras that were obtained recently. 
We investigate 
ascending chain conditions both for ideals and for 
Baxter ideals. The free Baxter algebras under 
consideration include free Baxter algebras on sets 
and free Baxter algebras on algebras. We also 
consider complete free Baxter algebras. 
\end{abstract}

\setcounter{section}{0}

\section{Introduction}

Let $C$ be a commutative ring and let $\lambda$ be an 
element of $C$. A {\bf Baxter $C$-algebra 
of weight $\lambda$} is a commutative $C$-algebra $R$ with a 
$C$-linear operator $P$ that satisfies the Baxter identity 
\margin{eq:bax}
\begin{equation}
 P(x)P(y)=P(xP(y))+P(yP(x))+\lambda P(xy), \forall x, y\in R. 
\label{eq:bax}
\end{equation}
The study of Baxter algebras was started by Baxter in 1963~\cite{Ba}.  
He was motivated by problems from fluctuation theory. 
In 1968, Rota~\cite{Ro1} began a systematic study of Baxter algebras 
from an algebraic point of view. 
Since then Baxter algebras have been related to hypergeometric 
functions, combinatorics, statistics, incidence algebras 
and theory of symmetric functions~\cite{Ro2,Ro3}.

Free Baxter algebras 
play a fundamental role in the study of Baxter algebras. 
Explicit descriptions of free Baxter algebras were 
first considered by Rota~\cite{Ro1} and Cartier~\cite{Ca}. 
In two recent papers~\cite{G-K1,G-K2}, 
William Keigher and the author 
furthered the work of Cartier and Rota, giving 
the explicit descriptions in complete generality. 
Using these constructions, further properties of Baxter algebras, 
in particular the zero divisors, were studied~\cite{Gu}, 
Baxter algebras were related to Hopf algebras~\cite{AGKO} 
and were applied to the umbral calculus~\cite{Gu2}. 

In this paper, we study ascending chain conditions in 
free Baxter algebras. 
Other than considering the noetherian ring property, 
we also consider modified noetherian properties, 
such as the the ascending chain condition for Baxter ideals. 
Let $X$ be a set. 
Denote $F_C(X)$ for the free Baxter $C$-algebra on $X$. 
The following is a summary of the main results on $F_C(X)$. 

\begin{enumerate}
\item
$F_C(\phi)$ is a noetherian ring if and only if 
$C$ is a noetherian $\QQ$-algebra (Theorem~\ref{thm:acc1}).
\item
If $X$ is not the empty set, then $F_C(X)$ is not noetherian
 (Theorem~\ref{thm:nacc1}). 
\item
If $C$ is a noetherian ring, then
$F_C(\phi)$ satisfies the ascending chain condition
for Baxter ideals (Theorem~\ref{thm:acc2}).
\item
If $X$ is not empty, 
then $F_C(X)$ 
of weight 0 does not satisfy the ascending chain condition for 
Baxter ideals. 
If $X$ is infinite, then $F_C(X)$ of any weight does not satisfy 
the ascending chain condition for Baxter ideals 
 (Corollary~\ref{co:nacc2}).
\end{enumerate}

As a generalization of free Baxter algebras on sets 
that were studied in \cite{Ro1} and \cite{Ca}, 
free Baxter algebras on $C$-algebras were introduced 
in \cite{G-K1}  
(see \S~\ref{sec:back} for more details). 
Ascending chain conditions in free Baxter algebras on 
$C$-algebras are also studied in this paper. 

In~\cite{G-K2}, we showed how one could complete a free Baxter 
algebra and get a complete free Baxter algebra. 
A summary of this construction is given in 
\S~\ref{sec:back}. 
This construction is similar to completing a free $C$-algebra 
(i.e., a polynomial ring with coefficients in $C$) and obtain a complete 
$C$-algebra (i.e, a power series ring with coefficients in $C$). 
In the current paper we also consider the ascending chain conditions 
in a complete free Baxter algebra.

We will provide some background on Baxter algebras in 
\S~\ref{sec:back}. 
In \S~\ref{sec:racc}, the ascending chain condition for ideals will
be studied and Theorems~\ref{thm:acc1} and \ref{thm:nacc1} will be 
proved. 
The ascending chain condition for Baxter ideals will be studied in 
\S~\ref{sec:bacc}. Theorem~\ref{thm:acc2}, 
Theorem~\ref{thm:nacc2} and Corollary~\ref{co:nacc2} 
are the main results in this section.

\section{Notations and background}
\margin{sec:back}
\label{sec:back}
We review concepts and results on Baxter algebras
that will be needed later in this paper.
See \cite{G-K1,G-K2,Gu} for detail.

\subsection{General notations}
We write $\NN$ for the set of natural numbers and 
$\NN_+=\{ n\in \NN\mid n>0\}$ for the positive integers. 

In this paper, every ring $C$ is commutative with
identity element 
$\bfone_C$, 
and every ring homomorphism
preserves the identity elements. 
For any $C$-modules $M$ and $N$,
the tensor product $M\otimes N$ is taken over $C$
unless otherwise indicated. 
For $n\in \NN$, denote the tensor power 
$\underbrace{M\otimes\ldots\otimes M}
    _{n\ {\rm factors}}$ by 
$M^{\otimes n}$
with the convention that $M^{\otimes 0}=C$.
This applies in particular if $M$ is a $C$-algebra. 
Let $\bfone$ be the identity element in a 
$C$-algebra $A$. 
We also use the notation 
$\bfone^{\otimes n} = \underbrace{\bfone \otimes\ldots\otimes \bfone}
    _{n\ {\rm factors}}$.

\subsection{Free Baxter algebras}
Let $(R,P)$ be a Baxter $C$-algebra of weight $\lambda$ 
with Baxter operator $P$. 
So $P$ satisfies the identity (\ref{eq:bax}). 
A {\bf Baxter ideal} of $(R,P)$ is an ideal $I$ of $R$ such
that $P(I)\subseteq I$.
The concepts of sub-Baxter algebras, quotient Baxter algebras
and homomorphisms of Baxter algebras can be easily defined. 

Let $A$ be a $C$-algebra. A {\bf free Baxter algebra} on $A$ 
is a Baxter algebra $(F_C(A),P_A)$ with a $C$-algebra homomorphism 
$j_A: A\to F_C(A)$ that satisfies the following universal property. 
For any Baxter $C$-algebra $(R,P)$ and any $C$-algebra
homomorphism $\varphi:A\rar R$, there exists a unique Baxter
$C$-algebra homomorphism $\tilde{\varphi}:(F_C(A),P_A)\rar
(R,P)$ such that the  diagram
\[ \begin{array}{ccc}
    A &\ola{j_A} & F_C(A)\\
    & {}_\varphi\!\!\searrow& \dap{\tilde{\varphi}}\\
    && R \end{array} \]
commutes.
Let $X$ be a set and let $A=C[X]$. 
Then $F_C(A)$ is the free Baxter algebra on $X$ in the usual 
sense. 
The existence of free Baxter algebras follows from the
general theory of universal algebras.
In order to get a good understanding of free Baxter algebras
and Baxter algebras in general,
it is desirable to find more explicit descriptions of
free Baxter algebras.

\subsection{Shuffle Baxter algebras}
\margin{sec:shuf}
\label{sec:shuf}
Motivated by the shuffle product of
iterated integrals~\cite{Re}, 
an explicit description of free Baxter algebras was given 
in~\cite{G-K1}. 
This generalizes earlier construction of 
free Baxter algebras by Cartier~\cite{Ca}. 
The resulting free Baxter algebras are called 
{\bf shuffle Baxter algebras}. 
We summarize the construction. 

For $m,n\in \NN_+$,
define the set of {\bf $(m,n)$-shuffles} by
\[ S(m,n)=
 \left \{ \sigma\in S_{m+n}
    \begin{array}{ll} {} \\ {} \end{array} \right .
\left |
\begin{array}{l}
\sigma^{-1}(1)<\sigma^{-1}(2)<\ldots<\sigma^{-1}(m),\\
\sigma^{-1}(m+1)<\sigma^{-1}(m+2)<\ldots<\sigma^{-1}(m+n)
\end{array}
\right \}.\]
Given an $(m,n)$-shuffle $\sigma\in S(m,n)$,
a pair of indices $(k, k+1)$,\ $1\leq k< m+n$, is
called an {\bf admissible pair} for $\sigma$
if $\sigma(k)\leq m<\sigma(k+1)$.
Denote $\calt^\sigma$ for the set of admissible pairs for $\sigma$.
For a subset $T$ of $\calt^\sigma$, we call the pair
$(\sigma,T)$ a {\bf mixable $(m,n)$-shuffle}.
Let $\mid T\mid$ be the cardinality of $T$.
We will identify $(\sigma,T)$ with $\sigma$
if $T$ is the empty set.
Denote
\[ \bs (m,n)=\{ (\sigma,T)\mid \sigma\in S(m,n),\
    T\subset \calt^\sigma\} \]
for the set of {\bf $(m,n)$-mixable shuffles}.
%Also denote
%\[ s(m,n)=\mid \bs(m,n)\mid .\]

Let $A$ be a $C$-algebra. 
For 
$x=x_1\otimes\ldots\otimes x_m\in A^{\otimes m}$,
$y=y_1\otimes \ldots\otimes y_n\in A^{\otimes n}$
and $(\sigma,T)\in \bar{S}(m,n)$,
the element
\[  \sigma (x\otimes y) =u_{\sigma(1)}\otimes u_{\sigma(2)} \otimes
    \ldots \otimes u_{\sigma(m+n)}\in A^{\otimes (m+n)},\]
where
\[ u_k=\left \{ \begin{array}{ll}
    x_k,& 1\leq k\leq m,\\
    y_{k-m}, & m+1\leq k\leq m+n, \end{array}
    \right . \]
is called a {\bf shuffle} of $x$ and $y$;
the element
\[ \sigma(x\otimes y; T)= u_{\sigma(1)}\hts u_{\sigma(2)} \hts
    \ldots \hts u_{\sigma(m+n)} \in A^{\otimes (m+n-\mid T\mid)},\]
where for each pair $(k,k+1)$, $1\leq k< m+n$,
\[ u_{\sigma(k)}\hts u_{\sigma(k+1)} =\left \{\begin{array}{ll}
    u_{\sigma(k)} u_{\sigma(k+1)},  &
     (k,k+1)\in T\\
    u_{\sigma(k)}\otimes u_{\sigma(k+1)}, &
    (k,k+1) \not\in T,
    \end{array} \right . \]
is called a {\bf mixable shuffle} of $x$ and $y$.

Fix a $\lambda\in C$ and a $C$-algebra $A$. 
There is a Baxter $C$-algebra of weight $\lambda$~\cite{G-K1}  

\[ \sha_C(A)=\sha_{C,\lambda}(A)= \bigoplus_{k\in\NN}
    A^{\otimes (k+1)}
= A\oplus A^{\otimes 2}\oplus \ldots \]
in which

\begin{itemize}
\item
the $C$-module structure
is the natural one,
\item
the multiplication is the {\bf mixed shuffle product},
defined by

\margin{eq:shuf}
\begin{equation}
x\shpr y\
=\sum_{(\sigma,T)\in \bs (m,n)} \lambda^{\mid T\mid }
    x_0y_0\otimes
\sigma(x^+\otimes y^+;T)
 \in \bigoplus_{k\leq m+n+1} A^{\otimes k}
\label{eq:shuf}
\end{equation}
for $x=x_0\otimes x_1\otimes\ldots \otimes x_m\in
A^{\otimes (m+1)}$ and
$y=y_0\otimes y_1\otimes\ldots\otimes y_n\in A^{\otimes (m+1)}$,
where $x^+=x_1\otimes \ldots \otimes x_m$ and
$y^+=y_1\otimes \ldots \otimes y_n$,
\item
the Baxter operator $P_A$ on
$\sha_C(A)$ is obtained by assigning
\[ P_A( x_0\otimes x_1\otimes \ldots \otimes x_n)
=\bfone_A\otimes x_0\otimes x_1\otimes \ldots\otimes x_n, \]
for all
$x_0\otimes x_1\otimes \ldots\otimes x_n\in A^{\otimes (n+1)}$.
\end{itemize}
$(\sha_C(A),P_A)$ is called the {\bf shuffle Baxter $C$-algebra on
$A$ of weight $\lambda$}.

For a given set $X$, we also let $(\sha_C(X),P_X)$ denote the
shuffle Baxter $C$-algebra $(\sha_C(C[X]),P_{C[X]})$, called the
{\bf shuffle Baxter $C$-algebra on $X$ (of weight $\lambda$).} Let
$j_A:A\rar \sha_C(A)$ (resp. $j_X:X\to \sha_C(X)$) be the
canonical inclusion map.

\begin{theorem}
\margin{thm:shua} 
\label{thm:shua} 
\cite{Ca,G-K1}
The pair $(\sha_C(A),P_A)$, together with the natural embedding $j_A$, is a
free Baxter $C$-algebra on $A$ of weight $\lambda$. 
Similarly, the pair 
$(\sha_C(X),P_X)$, together with the natural embedding
$j_X$,  is a free Baxter $C$-algebra on $X$
of weight $\lambda$.
\end{theorem}

We will use the following conventions in the rest of this paper. 
\begin{remark}
\begin{enumerate}
\item
Because of this theorem, we will use $\sha_C(X)$ instead of
$F_C(X)$ to denote the free Baxter $C$-algebra on $X$. 
\item
From the definition of the mixed shuffle product, we have 
\[ 
\!\!\!\!\!\!\!\! x\shpr y = \left \{ \begin{array}{ll}
    x_0 y_0, &\textrm{ if } x_0, y_0\in A,\\
x_0(y_0 \otimes y_1\otimes\ldots\otimes y_n), & \textrm{ if } 
   x=x_0\in A, y=y_0\otimes\ldots\otimes y_n\in A^{\otimes (n+1)}, 
\\
(x_0\otimes x_1\otimes\ldots\otimes x_n)y_0, & \textrm{ if } 
    x=x_0\otimes\ldots\otimes x_m\in A^{\otimes (m+1)}, 
    y=y_0\in A.
    \end{array} \right .
\]
This shows that the mixed shuffle product is compatible 
with the product in $A$.  
Thus we will suppress the symbol $\shpr$ in the mixed shuffle 
product unless there is the risk of confusion.  
\item
Unless otherwise specified, we use $A^{\otimes k}$ 
to denote the $C$-submodule of $\sha_C(A)$ instead 
of the tensor product algebra. 
\item
For $k\in \NN$, we denote  
$\Fil^k \sha_C(A)$ for $\bigoplus_{n\geq k} A^{\otimes (n+1)}$. 
\end{enumerate}
\end{remark}

\margin{ss:comp}
\subsection{Complete shuffle Baxter algebras}
\label{ss:comp}
We now take the completion of $\sha_C(A)$ in a manner similar to 
taking the completion of a polynomial ring to get a power series ring.

Given $k\in \NN_+$,
$\Fil^k \sha_C(A)$ 
is a Baxter ideal of $\sha_C(A)$.
On the other hand, consider the infinite product of $C$-modules
$\prod_{k\in \NN} A^{\otimes (k+1)}$.
It contains $\sha_C(A)$ as a dense subset with respect to the
topology defined by the filtration $\Fil^k \sha_C(A)$, $k\geq 0$.
All operations of the Baxter $C$-algebra $\sha_C(A)$ are continuous
with respect to this topology, hence extend uniquely to
operations on $\prod_{k\in\NN} A^{\otimes (k+1)}$,
making $\prod_{k\in\NN} A^{\otimes (k+1)}$ into a Baxter algebra
of weight $\lambda$.
We denote this Baxter algebra by 
$\widehat{\sha}_C(A)$ and denote the Baxter operator by $\hat{P}$. 
The pair $(\widehat{\sha}_C(A),\hat{P})$ is called the 
{\bf complete shuffle Baxter algebra on $A$}. 
It has been shown that $\widehat{\sha}_C(A)$ 
is a free object in the category of Baxter algebras
that are complete with respect to a canonical filtration
defined by the Baxter operator~\cite{G-K2}.

When $A=C$, we have
\[ \sha_C(C) = \bigoplus_{n\in \NN} C \bfone^{\otimes (n+1)}, 
\ \ \
 \widehat{\sha}_C(C)= \prod_{n\in\NN} C\bfone^{\otimes (n+1)},  
\]
where
$\bfone^{\otimes (n+1)}
= \underbrace{\bfone_C \otimes \ldots \otimes \bfone_C}
_{(n+1)-{\rm factors}}$. In this case the mixable shuffle
product formula~(\ref{eq:shuf}) gives

\margin{eq:unit}
\begin{equation}
\bfone^{\otimes (m+1)} \shpr \bfone^{\otimes (n+1)} =
\sum_{k=0}^m \binc{m+n-k}{n}\binc{n}{k} \lambda^k
\bfone^{\otimes (m+n+1-k)},\ \forall\ m,\ n\in \NN.
\label{eq:unit}
\end{equation}
This holds in both $\sha_C(C)$ and $\widehat{\sha}_C(C)$. 

\subsection{The internal construction}
\margin{sec:rota} 
\label{sec:rota} 
Later in the paper will use another 
construction of free Baxter algebras~\cite{G-K2},
generalizing the work of Rota~\cite{Ro1}. 
Since we will only need this construction in the special 
case when $A=C$, 
we will give a simplified description here. 
See~\cite{G-K2,Gu} for details. 

Define $ \frakA(C) = \prod_{n=1}^\infty C $
with componentwise addition and multiplication. 
Then $\frakA(C)$ is a $C$-algebra. It is in fact a 
Baxter $C$-algebra. 

\begin{prop}
\margin{prop:src}
\label{prop:src}
\cite{Gu}
Let $\lambda\in C$ be a non-zero divisor. 
Define 
$\Phi: \sha_C(C) \to \frakA(C)$ by sending 
$ b= \sum_{m=0}^\infty b_m \bfone^{\otimes (m+1)}\in \sha_C(C)$  to 
$\left ( \sum_{i=0}^{n-1}
    \bincc{n-1}{i}\lambda^i b_i \right)_{n\in\NN_+}\in \frakA(C)$.
Then $\Phi$ is an injective $C$-algebra homomorphism. 
Further, $\Phi$ extends to an injective $C$-algebra 
homomorphism 
$\Phi: \widehat{\sha}_C(C) \to \frakA(C)$. 
\end{prop}

\section{Ascending chain condition for ideals}
\margin{sec:racc}
\label{sec:racc}
In this section we prove the two theorems (Theorem~\ref{thm:acc1} 
and Theorem~\ref{thm:nacc1}) on the ascending chain condition 
for ideals in a free Baxter algebra. 

\subsection{The case when $A=C$}

\begin{theorem}
\margin{thm:acc1}
\label{thm:acc1}
\begin{enumerate}
\item
If $C$ is a noetherian $\QQ$-algebra, 
then $\sha_C(C)$ is a noetherian ring for every $\lambda\in C$.
\item
If $C$ is a noetherian $\QQ$-algebra and if $\lambda =0$,
then $\widehat{\sha}_C(C)$ is a noetherian ring.
\item
If $C$ is a $\QQ$-algebra,
$\lambda\in C$ is a not a zero divisor and
$\cap_{n\in\NN}\lambda^n C \neq 0$,
then $\widehat{\sha}_C(C)$ is not a noetherian ring. 
\item
If $C$ is not a $\QQ$-algebra, then
$\sha_C(C)$ and $\widehat{\sha}_C(C)$ are not noetherian rings.
\end{enumerate}
\end{theorem}

\proof
(1).
Let $C$ be a $\QQ$-algebra.
It is well-known that $R$ is a noetherian 
ring if and only if every ideal $I$ of $R$ is
finitely generated. 
So we only need to prove that any ideal of $\sha_C(C)$
is finitely generated.
The idea of the proof is the same as that of the
Hilbert basis theorem for $C[x]$. 
Let $I\subseteq \sha_C(C)$ be an ideal.
For each $j\in\NN$, let
\[\Sigma_j =\{ b_j\in C\mid \exists f_j\in I\ {\rm such\ that\ }
    f_j=\sum_{k=0}^j b_k \bfone^{\otimes (k+1)}\}. \]
Then for any $b_j,\ c_j\in \Sigma_j$, 
there are $f_j$ and $g_j$ in $I$ such that 
$f_j=\sum_{k=0}^j b_k \bfone^{\otimes (k+1)}$
and 
$g_j=\sum_{k=0}^j c_k \bfone^{\otimes (k+1)}$. 
So $f_j-g_j=\sum_{k=0}^j (b_k-c_k) \bfone^{\otimes (k+1)}$ 
is in $I$. Thus $b_j-c_j$ is in $\Sigma_j$. 
Also, for any 
$b_j\in \Sigma_j$ and $c\in C$, 
there is $f_j$ in $I$ such that $f_j=\sum_{k=0}^j b_k 
\bfone^{\otimes (k+1)}$. 
So $c f_j=\sum_{k=0}^j c b_k \bfone^{\otimes (k+1)}$ is in $I$. 
Thus $c b_j$ is in $\Sigma_j$. 
Therefore $\Sigma_j$ is an ideal of $C$.
Further, 
$ b_j\in \Sigma_j $ implies that 
there exists $f_j$ in $I$ such that 
$f_j=\sum_{k=0}^j b_k \bfone^{\otimes (k+1)}$. 
Thus $\bfone^{\otimes 2} f_j$ is in $I$. 
By equation~(\ref{eq:unit}), 
\begin{eqnarray*}
\lefteqn{\bfone^{\otimes 2} f_j = 
\sum_{k=0}^j b_k \bfone^{\otimes 2} \bfone^{\otimes (k+1)}} \\
&=& \sum_{k=0}^j b_k ((k+1)\bfone^{\otimes (k+2)}+k\bfone^{\otimes (k+1)})\\
&=& b_j (j+1) \bfone^{\otimes (j+2)} + {\rm\ lower\ degree\ terms }. 
\end{eqnarray*}
Here we define $\deg f=n$ if 
$f=\sum_{i=0}^\infty c_i \bfone^{\otimes (i+1)}$
with $c_n\neq 0$ and $c_i=0$ for $i>n$ and define $\deg 0=\infty$.
Thus $b_j (j+1)$ is in $\Sigma_{j+1}$.
Since $\Sigma_{j+1}$ is an ideal and $C$ is a $\QQ$-algebra,
we have $b_j= (j+1)^{-1} b_j (j+1)\in \Sigma_{j+1}$. 
Thus $\Sigma_j\subseteq \Sigma_{j+1}$.
Since $C$ is noetherian, this chain of ideals stabilizes,
say at $j=m$. Then $\Sigma_m=\cup_{j=1}^\infty \Sigma_j$,
and is finitely generated. 
Let $b_1^{(m)},\ldots,b_{k_m}^{(m)}$ be a set of generators of
$\Sigma_m$. Then we have $f_i^{(m)}\in I,\ 1\leq i\leq k_m$, with 
$f_i^{(m)}=b_i^{(m)} \bfone^{\otimes (m+1)} +g_i^{(m)},\
\deg g_i^{(m)}<m$. 
For each $j<m$, $\Sigma_j$ is also finitely generated with
a set of generators $b_1^{(j)},\ldots,b_{k_j}^{(j)}$.
Then there are $f_i^{(j)}\in I$ such that
$f_i^{(j)}= b_i^{(j)}\bfone^{\otimes (j+1)}+g_i^{(j)}\in I$
with $\deg g_i^{(j)}<j$.
To prove the theorem, we only need to prove that 
$I$ is the ideal generated by

\[\{f_1^{(0)},\ldots,f_{k_0}^{(0)},f_1^{(1)},\ldots,f_{k_1}^{(1)},
\ldots, f_1^{(m)},\ldots,f_{k_m}^{(m)}\}.\]
Let $I\tp$ be the ideal generated by this set.
Clearly $0$ is in $I\tp$. For $f\in I$ with $f\neq 0$, 
we use induction on $\deg f$ to show that $f$ is in $I\tp$. 
If $f\in I$ with $\deg f=0$, then $f\in \Sigma_0 \bfone=\Sigma_0$,
so can be expressed as a $C$-linear combination of
$f_1^{(0)},\ldots,f_{k_0}^{(0)}$. Thus $f\in I\tp$.
Now for any $n>0$. 
Assume that all $f\in I$ with $\deg f<n$ are in $I\tp$ and 
take $f\in I$ with $\deg f=n$. Write
$f=b_n \bfone^{\otimes (n+1)}+g,\ b_n\neq 0,\ \deg g<n$.
Then $b_n$ is in $\Sigma_n$. 
If $n\geq m$, then by the definition of $m$,
$\Sigma_n=\Sigma_m$. So
$b_n=\sum_i a_i b_i^{(m)}$ for some $a_i\in C$. 
Consider 
$h=\bfone^{\otimes (n-m+1)} \sum_i a_i f_i^{(m)}$. 
From $f_i^{(m)}\in I$ we see that $h$ is in $I\tp$. 
Also, by equation~(\ref{eq:unit})
\begin{eqnarray*}
\lefteqn{ \bfone^{\otimes (n-m+1)} f_i^{(m)} 
= b_i^{(m)} \bfone^{\otimes (n-m+1)} \bfone^{\otimes (m+1)} 
 + \bfone^{\otimes (n-m+1)} g^{(m)}_i }\\
&=& b_i^{(m)} \binc{n}{m} \bfone^{\otimes (n+1)} + 
{\rm\ lower\ degree\ terms }. 
\end{eqnarray*}
Thus 
\[ h= \sum_i a_i \bfone^{\otimes (n-m+1)} f_i^{(m)} 
= \binc{n}{m} \sum_i a_i b_i^{(m)} \bfone^{\otimes (n+1)} 
+ {\rm\ lower\ degree\ terms }. \]
Therefore $\binc{n}{m}^{-1} h$, still in $I\tp$ since $C$ is a 
$\QQ$-algebra, 
has the same leading coefficient as $f$. 
Hence 
$f- \binc{n}{m}^{-1} h$ 
has degree less then $n$.
Since $f- \binc{n}{m}^{-1} h $ is in $I$, it is in
$I\tp$ by induction. Then $f$ is in $I\tp$.
If $n<m$, then $b_n=\sum_i a_i b_i^{(n)},\ a_i\in C$. 
So $\sum_i a_i f_i^{(n)}$ is in $I\tp$ and has the same
leading coefficient as $f$.
Then $f-\sum_i a_i f_i^{(n)}\in I$ with
$\deg(f-\sum_i a_i f_i^{(n)})<n$.
By induction, $f-\sum_i a_i f_i^{(n)}$ is in $I\tp$.
Hence $f$ is in $I\tp$.

\mvp
\noindent
(2)
We now consider $\widehat{\sha}_C(C)$.
Since $\lambda=0$, equation~(\ref{eq:unit}) becomes 
\[ \bfone^{\otimes (m+1)} \bfone^{\otimes (n+1)} 
    = \binc{m+n}{n} \bfone^{m+n+1}. \]
Then we can use an argument that is similar to the previous part of the 
proof. Just replace 
$\Sigma_j$ by
\[ \Omega_j=\{b_j\in C\mid \exists f_j\in I {\rm\ such\ that\ }
    f_j=\sum_{k=j}^\infty  b_k \bfone^{\otimes (k+1)} \} \]
and follow the well-known argument in proving 
that $C[[x]]$ is a noetherian ring.

\mvp
\noindent
(3)
We only need to find ideals $I_n$, $n\geq 1$, of 
$\widehat{\sha}_C(C)$ such that, for each $n$, 
$I_n\subsetneq I_{n+1}$. For this purpose, we will construct a 
sequence $d^{(k)},\ k\geq 1$ of elements in $\widehat{\sha}_C(C)$ 
with the property that, for each $m\geq 1$, 
\margin{eq:ann}
\begin{equation}
d^{(k)} \Fil^{m}\widehat{\sha}_C(C) \left \{ 
\begin{array}{ll} \neq 0, & {\rm\ if\ } m=k-1, \\
= 0, & {\rm\ if\ } m=k. \end{array} \right . 
\label{eq:ann}
\end{equation}
We then let $I_n$ be the ideal of $\hatsha_C(C)$ generated by 
$d^{(k)},\ 1\leq k\leq n$. Then clearly $I_n\subseteq I_{n+1}$. 
From equation~(\ref{eq:ann}) we have 
\[
I_n \Fil^{m}\widehat{\sha}_C(C) \left \{ 
\begin{array}{ll} \neq 0, & {\rm\ if\ } m=n-1, \\
= 0, & {\rm\ if\ } m=n. \end{array} \right . 
\]
In particular, $I_n \Fil^n\hatsha_C(C)= 0$ 
while $I_{n+1} \Fil^n\hatsha_C(C)\neq 0$. 
Therefore, $I_n\neq I_{n+1}$, as is desired. 
The rest of the proof will be devoted to the construction 
of such a sequence $d^{(k)},\ k\geq 1$. 

We first assume that $C$ is $\QQ(x)=\QQ[x,x^{-1}]$ 
and assume that 
the weight of 
$\hatsha_{\QQ(x)}(\QQ(x))$ is $x$.
%Then $C=S^{-1} \QQ[x]$, the localization of the polynomial
%algebra $\QQ[x]$ at the multiplicative set $S=\{x^n, n\geq 0\}$.
For a fixed $k\in \NN_+$, we want to find a solution 
$b=b^{(k)}\in \hatsha_{\QQ(x)}(\QQ(x))$ to the equation
\margin{eq:zero1}
\begin{equation}
 \bfone^{\otimes (k+1)} b =0. 
\label{eq:zero1}
\end{equation}
Write $b= \sum_{n=0}^\infty b_n \bfone^{\otimes (n+1)},\ b_n\in 
\QQ(x)$.
We have
\begin{eqnarray*}
\lefteqn{ \bfone^{\otimes (k+1)} b 
= \bfone^{\otimes (k+1)}
    (\sum_{n=0}^\infty b_n \bfone^{\otimes (n+1)})}\\
&=& \sum_{n=0}^\infty b_n \bfone^{\otimes (k+1)}\bfone^{\otimes (n+1)}\\
&=& \sum_{n=0}^\infty b_n
\left (\sum_{i=0}^k \bincc{n+k-i}{k}\bincc{k}{i} x^i \bfone^{\otimes (n+k-i+1)}
    \right )
\ \ \ ({\rm equation~(\ref{eq:unit})})\\
&=&\!\!\!\! \sum_{i=0}^k \sum_{n=0}^\infty \bincc{n+k-i}{k}\bincc{k}{i}
    x^i b_n \bfone^{\otimes (m+k-i+1)}
\ ({\rm exchanging\ the\ order\ of\ summation}) \\
&=& \sum_{i=0}^k \sum_{m=k-i}^\infty \bincc{m}{k}\bincc{k}{i}
    x^i b_{m-k+i} \bfone^{\otimes (m+1)} 
\ ({\rm\ replacing\ } n {\rm\ by\ } m-k+i)\\
&=& \sum_{i=0}^k \sum_{m=k}^\infty \bincc{m}{k}\bincc{k}{i}
    x^i b_{m-k+i} \bfone^{\otimes (m+1)} 
\ \ \ \left (\binc{m}{k}=0 {\rm\ for\ } m<k \right)\\
&=&\!\!\!\! \sum_{m=k}^\infty \bincc{m}{k}
\left (\sum_{i=0}^k \bincc{k}{i} x^i b_{m-k+i} \right)
    \bfone^{\otimes (m+1)} 
\  ({\rm  exchanging\ the\ order\ of\ summation}).
\end{eqnarray*}
Thus finding a solution $b\in \hatsha_{\QQ(x)}(\QQ(x))$ of 
equation~(\ref{eq:zero1}) is equivalent to finding 
solutions $b_n\in \QQ(x)$ of the system of equations 
\margin{eq:zero2'}
\begin{equation}
\bincc{m}{k} \sum_{i=0}^k \bincc{k}{i} x^i b_{m-k+i} =0,\ m\geq k. 
\label{eq:zero2'}
\end{equation}
Since $\QQ(x)$ has characteristic zero, solving system~(\ref{eq:zero1})
is equivalent to solving the system 
\margin{eq:zero2}
\begin{equation}
\sum_{i=0}^k \bincc{k}{i} x^i b_{m-k+i} =0,\ m\geq k
\label{eq:zero2}
\end{equation}
in $\QQ(x)$. 
This last system of equations can be rewritten as 
\margin{eq:zero3}
\begin{equation}
b_m= -x^{-k} \left (\sum_{i=0}^{k-1}
    \bincc{k}{i}x^i b_{m-k+i} \right ),\ m\geq k.
\label{eq:zero3}
\end{equation}
For $m=k$, we have
\[b_k= -x^{-k} \left (\sum_{i=0}^{k-1}
    \bincc{k}{i}x^i b_{i} \right ),\ m\geq k.\]
Choosing $b_0=1$ and $b_i=0,\ 1\leq i\leq k-1$, we have
$b_k=b_k(x) = -x^{-k}$.
Inductively, these values of $b_0,\ldots, b_{k-1}$ and 
equation~(\ref{eq:zero3}) uniquely determine
a rational function $b^{(k)}_m(x)\in \QQ(x)$ for each $m\geq k$, 
giving us a non-zero solution 
\[ \{b_m=b^{(k)}_m(x),\ m\geq 0\} \]
of the linear system~(\ref{eq:zero2}) with values in $\QQ(x)$. 
Hence we obtain a non-zero solution

\[ b^{(k)}=\sum_{n=0}^\infty b^{(k)}_n (x) \bfone^{\otimes (n+1)}  \]
to equation~(\ref{eq:zero1}) in $\hatsha_{\QQ(x)}(\QQ(x))$.
Note that $b^{(k)}$ is a function of $x$. We denote it 
by $b^{(k)}(x)$. 

Using Proposition~\ref{prop:src}, from 
$\bfone^{\otimes (k+1)} b^{(k)}(x)=0$ we have
\margin{eq:xzero0}
\begin{equation}
\Phi(\bfone^{\otimes (k+1)})\Phi(b^{(k)}(x))=
\Phi(\bfone^{\otimes (k+1)} b^{(k)}(x))=0
\label{eq:xzero0}
\end{equation}
in $\frakA(C)$.
But from the definition of $\Phi$, we have
\[ \Phi\left (\bfone^{\otimes (k+1)}\right ) =\left (\bincc{n-1}{k} x^k\right )_n\]
and
\[ \Phi\left (b^{(k)}(x)\right)= \left (\sum_{i=0}^{n-1} b^{(k)}_i(x) \bincc{n-1}{i} x^i\right )_n. \]
So

\[ \Phi(\bfone^{\otimes (k+1)})\Phi(b^{(k)}(x))
= \left (\binc{n-1}{k} x^k \left (\sum_{i=0}^{n-1} b^{(k)}_i(x) \bincc{n-1}{i} 
x^i\right )
\right )_n,
\]
and equation~(\ref{eq:xzero0}) becomes  
\[ \left (\bincc{n-1}{k} x^k 
\left (\sum_{i=0}^{n-1} b^{(k)}_i(x) \bincc{n-1}{i} x^i\right)
\right )_n =0.\]
Therefore, 
\[ \bincc{n-1}{k} x^k 
\left (\sum_{i=0}^{n-1} b^{(k)}_i(x) \bincc{n-1}{i} x^i\right)
    =0,\ n\geq 1.\]
Since $\bincc{n-1}{k}\neq 0$ for $n\geq k+1$, we must have

\margin{eq:xzero}
\begin{equation}
\sum_{i=0}^{n-1} b^{(k)}_i(x) \bincc{n-1}{i} x^i =0,\  
n\geq k+1.
\label{eq:xzero}
\end{equation}

We now let $C$ be any $\QQ$-algebra and let $\lambda$ be a
non-zero divisor in $C$. Let $S=\{\lambda^n,\ n\geq 0\}$ 
and consider the localization $S^{-1}C$. 
Since $\lambda$ is not a zero divisor, the assignment 
$x\mapsto \lambda$ induces a ring homomorphism 
\[ C[x,x^{-1}] \to S^{-1} C. \] 
Let $b^{(k)}_i(\lambda)$ be the image of $b^{(k)}_i(x)$ under this 
homomorphism. 
Then from equation~(\ref{eq:xzero}) we have the equations 

\[\sum_{i=0}^{n-1} b^{(k)}_i(\lambda) \bincc{n-1}{i} \lambda^i =0,  
n\geq k+1\]
in $S^{-1}C$. 
This shows that, for

\[ b^{(k)}(\lambda)=
\sum_{n=0}^\infty b^{(k)}_n(\lambda) \bfone^{\otimes (n+1)},\]
the $n$-th component of
$\Phi(b^{(k)}(\lambda))$ is zero for $n\geq k+1$. 
On the other hand, from 
the definition of $\Phi$, 
for any $\alpha\in \Fil^{k} \widehat{\sha}_{S^{-1}C}(S^{-1}C)$,
the $n$-th component of
$\Phi (\alpha)$ is zero for $n\leq k$. 
Since the product in $\frakA(C)$ is defined componentwise, 
we further have, 
for $\alpha\in \Fil^{k} \widehat{\sha}_{S^{-1}C}(S^{-1}C)$,
\[
\Phi(\alpha b^{(k)}(\lambda))=\Phi(\alpha)\Phi(b^{(k)}(\lambda))=0.
\]
Since $\Phi$ is injective, we have 
\margin{eq;qzero}
\begin{equation}
b^{(k)}(\lambda) \Fil^{k}\widehat{\sha}_{S^{-1}C}(S^{-1}C)=0. 
\label{eq:qzero}
\end{equation}

By the assumption of the theorem, there is a non-zero element $c$ in
$\cap_{n=0}^\infty \lambda^n C$. 
Fix such a $c$ and define 
\[ d^{(k)}= c b^{(k)}(\lambda),\ k\geq 1.\]
Then $d^{(k)}$ is in $\widehat{\sha}_{S^{-1}C}(S^{-1}C)$. 
To finish the proof, we only need to show that each $d^{(k)}$ is  
in $\widehat{\sha}_C(C)$ and satisfies equation~(\ref{eq:ann}). 
Here we regard $\widehat{\sha}_C(C)$ 
as the subalgebra of 
$\widehat{\sha}_{S^{-1}C}(S^{-1}C)$ 
consisting of sequences 
$\sum_{n=0}^\infty a_n \bfone^{\otimes (n+1)}$ with 
$a_n\in C$, $n\geq 0$. 
This is justified 
because $C$ can be identified with a subalgebra of 
$S^{-1}C$ since $\lambda$ is not a zero divisor, 
and because 
$\widehat{\sha}_C(C)= \prod_{n=0}^\infty C \bfone^{\otimes (n+1)}$ 
and 
$\widehat{\sha}_{S^{-1}C}(S^{-1}C)= 
    \prod_{n=0}^\infty S^{-1}C \bfone^{\otimes (n+1)}$, 
as we have seen in \S~\ref{ss:comp}.

Since $c$ is in $\cap_{n=0}^\infty \lambda^n C$, 
for each $n\in \NN$, there is $c_n\in C$ such that
$c=\lambda^n c_n$.
Further, 
for each $n\geq 0$, the rational function 
$x^n b^{(k)}_n(x)$ in $\QQ(x)$ is a polynomial in $\QQ[x]$.
This is clear for $0\leq n\leq k$ and the general case 
follows by induction on $n$. 
Thus for each $n\geq 0$, 
$ cb^{(k)}_n(\lambda) = 
 c_n \lambda^n b^{(k)}_n(\lambda)$
is an element in $C$.
Therefore, 
\[ d^{(k)}=cb^{(k)}(\lambda) = 
\sum_{n=0}^\infty c_n \lambda^n b^{(k)}_n(\lambda) \bfone^{\otimes (n+1)}
    \]
is an element in $\widehat{\sha}_C(C)$.

Further,

\[ c b^{(k)}(\lambda) \Fil^{k}\widehat{\sha}_C(C)
=c (b^{(k)}(\lambda) \Fil^{k}\widehat{\sha}_C(C))
=0 \]
since we have proved that 
$b^{(k)}(\lambda) \Fil^{k}\widehat{\sha}_{S^{-1}C}(S^{-1}C)=0$ 
in equation~(\ref{eq:qzero}). 
On the other hand, since we have chosen 
$b^{(k)}_0 (x) =1$ and $b^{(k)}_1(x)=\ldots =b^{(k)}_{k-1}(x)=0$ 
in
$b^{(k)}=b^{(k)}(x)=\sum_{n=0}^\infty b^{(k)}_n (x) \bfone^{\otimes (n+1)}$
we see that, for $1\leq n\leq k$, 
the $n$-th component of $\Phi(b^{(k)})$ is 

\[ \sum_{i=0}^{n-1}b^{(k)}_i(x) \bincc{n-1}{i} x^i
= 1.\]
Thus, the $n$-th component of $\Phi(cb^{(k)}(\lambda))$ for
$1\leq n\leq k$ is $c$.
For $\bfone^{\otimes k}\in C^{\otimes k}$, 
the $k$-th component of $\Phi(\bfone^{\otimes k})\in \frakA(C)$ 
is $\lambda ^{k-1}$.
Thus the $k$-th component of
$\Phi(cb^{(k)}(\lambda) \bfone^{\otimes k})$ is
$c\lambda^{k-1}$.
It is not zero, since $c$ is not zero and $\lambda$ is not a zero divisor.
Thus $cb^{(k)}(\lambda) \Fil^{k-1} \widehat{\sha}_C(C)$ is not
zero.
Thus we have shown that the elements 
$d^{(k)}=c b^{(k)}(\lambda),\ k\geq 1$, of 
$\widehat{\sha}_C(C)$
satisfy equation~(\ref{eq:ann}). 
This completes the proof of part 3 of the theorem.

\vspace{0.3cm}
\noindent
(4) 
If $C$ is not a $\QQ$-algebra, 
then there is a prime number
$p$ such that $p\cdot \bfone_C$ is not a unit in $C$.
Thus there is a maximal ideal $M$ of $C$ containing $p\cdot\bfone_C$.
Let $F=C/M$ be the residue field. Then
$F$ is an algebra over the finite field $\FF_p$.
Let $\tilde{M}$ be the Baxter ideal of $\sha_C(C)$
generated by $M$. Then by Proposition~3.3 in~\cite{Gu}, 
\[\sha_C(C)/\tilde{M} \cong \sha_C(F) \cong \sha_F(F).\]
If $\sha_C(C)$ were noetherian, then its quotient
$\sha_F(F)$ would also be noetherian. 
Thus the theorem follows from the following lemma.
\proofend

\begin{lemma}
\margin{lem:fin}
\label{lem:fin}
If $F$ is a field of non-zero characteristic $p$,
then $\sha_F(F)$ is not a noetherian ring.
\end{lemma}

\proof
For each $k\geq 1$, define
\[ I_k = \sum_{n}{}\tp F\bfone^{\otimes (n+1)} 
    \subseteq \sha_F(F),\]
where the sum is over all $n\in \NN$ with $p^k\nmid n$. 
We prove that each $I_k$ is an ideal of $\sha_F(F)$.
For this we only need to show that
$\bfone^{\otimes (m+1)} 
\bfone^{\otimes (n+1)}\in I_n$ for $m\in\NN$ and $p^k\nmid n$. 
We have
\[ \bfone^{\otimes (m+1)} \bfone^{\otimes (n+1)}
    =\sum_{i=0}^n \binc{m+n-i}{n}\binc{n}{i}
    \bfone^{\otimes (m+n-i+1)}.\]
For each $0\leq i\leq n$, 
if $p^k\nmid m+n-i$, then $\bfone^{\otimes (m+n-i+1)}\in I_n$;
if $p^k\mid m+n-i$, then from $p^n\nmid n$ we have
$\binc{m+n-i}{n}\equiv 0\pmod{p}$ \cite[p.68]{Kn}. 
So $\bfone^{\otimes (m+1)} \bfone^{\otimes (n+1)}\in I_n$, 
and $I_n$ is an ideal.
By definition we have $\bfone^{\otimes (p^n+1)}\in I_{n+1}$ but
$\bfone^{\otimes (p^n+1)}\not\in I_n$ for each 
$n\geq 1$.
Therefore, $I_n$ is a strictly increasing sequence of
ideals, as needed.
\proofend

\subsection{The general case}

\begin{theorem}
\margin{thm:nacc1}
\label{thm:nacc1}
Let $C$ be a ring of characteristic zero. For any non-empty set 
$X$, the free Baxter algebra $\sha_C(X)$ is not a noetherian algebra. 
\end{theorem}

\proof
We start with the case when $X$ is a singleton $\{x\}$. 
For each integer $n\geq 1$, let $\Sigma_n$ be the ideal of 
$\sha_C(X)$ generated by $1\otimes x^i,\ 1\leq i\leq n$. 
To prove the theorem, it suffices to show that 
$\Sigma_n\subsetneq \Sigma_{n+1}$ for each $n\geq 1$. 
We prove this by contradiction. Assume that 
$\Sigma_{n+1} =\Sigma_n$ for some $n$. 
Then in particular, $1\otimes x^{n+1}\in \Sigma_n$. 
Thus $1\otimes x^{n+1}$ can be expressed in the form 
\[ \sum_{k=1}^n (1\otimes x^k) G_k,\ G_k\in \sha_C(X).\]
The construction of $\sha_C(X)$ shows that $\sha_C(X)$ is 
a free $C[x]$-module on the set 
\[ \calx=\{1\}\cup 
    \{  1\otimes x^{i_1}\otimes x^{i_2}\otimes \ldots 
\otimes x^{i_m}\mid 
i_j \in\NN {\rm \ for\ } 1\leq j\leq m,\ m\geq 1 \}. \]
Define 
\[ \cali = \{\phi\} \cup \left ( 
    \bigcup_{m=1}^\infty \NN^m \right )\]
and, for $I\in \cali$, denote 

\[ x^I= \left \{ \begin{array}{ll} 
    1, & {\rm\ if\ } I=\phi, \\
    1\otimes x^{i_1}\otimes x^{i_2}\otimes \ldots \otimes x^{i_m}, 
    & {\rm\ if\ } I=(i_1,\ldots,i_m). \end{array} 
    \right .\]
Then 
$ \calx =\{ x^I \mid I\in \cali \} $.
Thus each $G_k$ above, $1\leq k\leq n,$ can be written as 
$ \sum_{I\in \cali} g^{(k)}_I x^I$ for unique 
$g^{(k)}_I \in C[x]$ and we have 
\margin{eq:com1}
\begin{equation}
 1\otimes x^{n+1} 
= \sum_{k=1}^{n} (1\otimes x^k) (\sum_{I\in \cali} g^{(k)}_I x^I)
= \sum_{k=1}^{n} \sum_{I\in \cali} g^{(k)}_I (1\otimes x^k)\ x^I. 
\label{eq:com1}
\end{equation}
We will derive a contradiction from this equation. 

Since elements in $\calx$ form a basis for the free $C[x]$-module 
$\sha_C(X)$, we can write  
\margin{eq:com2}
\begin{equation} 
\sum_{k=1}^{n} \sum_{I\in \cali} g^{(k)}_I (1\otimes x^k)\ x^I
= \sum_{J\in \cali} h_J x^J
\label{eq:com2}
\end{equation}
for unique $h_J\in C[x],\ J\in \cali$. 
Comparing this with equation~(\ref{eq:com1}), we see that  
$h_J=1$ if $J=n+1$ and $h_J=0$ for all other $J\in \cali$. 
For $J\in \cali$, define 
\[ \mid J \mid = \left \{ \begin{array}{ll} 
 0, & {\rm\ if\ } J=\phi, \\
j_1+\ldots + j_m, & {\rm\ if\ } J=(j_1,\ldots,j_m). 
\end{array} \right .\]
Then we in particular have $h_J=0$ for 
$\mid J \mid \neq n+1$. 
Thus equation~(\ref{eq:com2}) becomes 
\margin{eq:com2'}
\begin{equation}
\sum_{k=1}^{n} \sum_{I\in \cali} g^{(k)}_I (1\otimes x^k)\ x^I
=\sum_{\mid J\mid =n+1} h_J x^J
\label{eq:com2'}
\end{equation}
and equation~(\ref{eq:com1}) becomes 

\margin{eq:com3}
\begin{equation}
1\otimes x^{n+1}=\sum_{\mid J\mid =n+1} h_J x^J.
\label{eq:com3}
\end{equation}
Next we will study the relation between the 
coefficients $g^{(k)}_I$ and $h_J$ more carefully. 

Fix a $k\in \NN$ and an $I\in \cali$. From the definition of 
the mixable shuffle product in equation~(\ref{eq:shuf}), 
we have 
\[ (1\otimes x^{k}) x^I = 1\otimes x^{k}\]
when $I=\phi$; while 
when $I=(i_1,\ldots,i_m)\in \cali$, we have 
\begin{eqnarray*}
 (1\otimes x^{k}) x^I
&=& 1\otimes ( x^{k}\otimes x^{i_1}\otimes \ldots \otimes 
x^{i_m} 
+ x^{i_1}\otimes x^{k} \otimes x^{i_2}\ldots \otimes x^{i_m} \\
&& + \ldots 
+ x^{i_1}\otimes \ldots x^{i_{m-1}}\otimes x^{k} \otimes x^{i_m} 
+ x^{i_1}\otimes \ldots \otimes x^{i_m}\otimes x^{k}
\\
&&  +
\lambda (x^{k+i_1}\otimes x^{i_2}\otimes \ldots \otimes x^{i_m} 
+ x^{i_1}\otimes x^{k+i_2}\otimes x^{i_3}\otimes\ldots\otimes 
    x^{i_m} \\
&&+ \ldots 
+ x^{i_1}\otimes \ldots \otimes x^{i_{m-1}}\otimes x^{k+i_m})).  
\end{eqnarray*}
Note that for each of the basis elements $x^J\in \calx$ 
that occurs 
on the right hand side of the equation, we have 
$\mid J\mid =k+\mid I\mid$. 
This shows that in equation~(\ref{eq:com2'}), 
a coefficient
$h_J$ on the right hand side must be a sum of the 
coefficients $g^{(k)}_I$ on the left hand side 
with the property $\mid J\mid =k+\mid I\mid$. 
Thus equation~(\ref{eq:com2'}) becomes 

\[\sum_{k=1}^{n}\  \sum_{\mid I\mid=n+1-k} 
    g^{(k)}_I (1\otimes x^{k})\ x^I
=\sum_{\mid J\mid =n+1} h_J x^J. \]
Exchanging the order of summation on the left hand side, 
we have 
\margin{eq:com3'}
\begin{equation}
\sum_{m=1}^{n} \ \sum_{\mid I\mid=m,\\ k=n+1-m} 
    g^{(k)}_I (1\otimes x^{k})\ x^I
=\sum_{\mid J\mid =n+1} h_J x^J
\label{eq:com3'}
\end{equation}
Since $\mid I\mid =n+1-k$ and 
$1\leq k\leq n$, we have 
$\mid I\mid\neq 0$ for any $I$ in this equation. 
So $I\neq \phi$ and hence 
$I=(i_1,\ldots,i_m)$ for $i_j\in \NN$ and $m\geq 1$. 
Then we have 

\begin{eqnarray*}
\lefteqn{\sum_{\mid I\mid =m, k=n+1-m} 
    g^{k}_I (1\otimes x^{k}) x^I }\\
&=& \sum_{k+i_1+\ldots+i_m=n+1} g^{(k)}_{i_1,\ldots,i_m} 
(1\otimes x^{k})\ 
    (1\otimes x^{i_1}\otimes x^{i_2}\otimes \ldots \otimes x^{i_m}) \\
&=&
\sum_{k+i_1+\ldots+i_m=n+1} g^{(k)}_{i_1,\ldots,i_m} 
\otimes \left (  x^{k}\otimes x^{i_1}\otimes \ldots \otimes 
x^{i_m} 
+ x^{i_1}\otimes x^{k} \otimes x^{i_2}\ldots \otimes x^{i_m} 
 \right . \\
&& + \ldots 
+ x^{i_1}\otimes \ldots x^{i_{m-1}}\otimes x^{k} \otimes x^{i_m} 
+ x^{i_1}\otimes \ldots \otimes x^{i_m}\otimes x^{k}
\\
&&  +
\lambda (x^{k+i_1}\otimes x^{i_2}\otimes \ldots \otimes x^{i_m} 
+ x^{i_1}\otimes x^{k+i_2}\otimes x^{i_3}\otimes\ldots\otimes 
    x^{i_m} \\
&&\left . + \ldots 
+ x^{i_1}\otimes \ldots \otimes x^{i_{m-1}}\otimes x^{k+i_m})
\right ). 
\end{eqnarray*}
Let $G_m^{(2)}$ (resp. $G_m^{(1)}$) 
be the sum of the terms on the right hand side  
in which the tensor product has $m+2$ (resp. $m+1$) tensor factors. 
More precisely, 
\begin{eqnarray*}
G_m^{(2)}
&= & 
\sum_{k+i_1+\ldots+i_m=n+1} g^{(k)}_{i_1,\ldots,i_m} 
\otimes  (  x^{k}\otimes x^{i_1}\otimes \ldots \otimes 
x^{i_m} 
+ x^{i_1}\otimes x^{k} \otimes x^{i_2}\ldots \otimes x^{i_m}\\
&& + \ldots 
+ x^{i_1}\otimes \ldots x^{i_{m-1}}\otimes x^{k} \otimes x^{i_m} 
+ x^{i_1}\otimes \ldots \otimes x^{i_m}\otimes x^{k}), \\
G_m^{(1)} &= & 
\sum_{k+i_1+\ldots+i_m=n+1} \lambda g^{(k)}_{i_1,\ldots,i_m} 
\otimes  
(x^{k+i_1}\otimes x^{i_2}\otimes \ldots \otimes x^{i_m} \\
&& + x^{i_1}\otimes x^{k+i_2}\otimes x^{i_3}\otimes\ldots\otimes 
    x^{i_m} 
 + \ldots 
+ x^{i_1}\otimes \ldots \otimes x^{i_{m-1}}\otimes x^{k+i_m}).
\end{eqnarray*}
Then from equation~(\ref{eq:com3}) and equation~(\ref{eq:com3'}), 
we have 

\margin{eq:com4}
\begin{equation}
1\otimes x^{n+1} 
=\sum_{m=1}^n (G_m^{(2)}+ G_m^{(1)}). 
\label{eq:com4}
\end{equation}
Thus for each $r\geq 2$, the sum of the terms on the right 
hand side of equation~(\ref{eq:com4}) 
with $r$ tensor factors is given by 
\[ \left \{ \begin{array}{ll}
0, & {\rm\ when\ } r< 2, \\
G^{(1)}_1, & {\rm\ when\ } r=2, \\
 G_{r-2}^{(2)}+G_{r-1}^{(1)}, & {\rm\ when\ } 3\leq r\leq n+1, \\
G^{(2)}_n, & {\rm\ when\ } r=n+2, \\
0, & {\rm\ when\ } r>n+2. 
\end{array} \right . \]
Therefore from equation~(\ref{eq:com4}) we have 
\margin{eq:com5}
\begin{equation}
\left \{ \begin{array}{ll}
G^{(1)}_1 = 1\otimes x^{n+1}, \\
G_{r-2}^{(2)}+G_{r-1}^{(1)} =0 , & 3\leq r\leq n+1,\\
G^{(2)}_n = 0. 
\end{array} \right . 
\label{eq:com5}
\end{equation}
From the definition of $G^{(2)}_m$, we see that the sum 
of the coefficients of all the basis elements 
$x^I\in \calx$ in $G^{(2)}_m$ is 
$ (m+1) g_m$ 
where 
\[ g_m= \sum_{k+i_1+\ldots+i_m=n+1} g^{(k)}_{i_1,\ldots,i_m} . \] 
Similarly, the sum of the coefficients of all the basis 
elements $x^I\in \calx$ in $G^{(1)}_m$ is $\lambda mg_m$. 
Therefore, the sum of the coefficients of all monomials in 
$G^{(2)}_{r-2} +G^{(1)}_{r-1}$ is

\margin{eq:com6}
\begin{equation}
\left \{ \begin{array}{ll}
\lambda g_1, & {\rm\ when\ } r=2, \\
(r-1)g_{r-2} +\lambda (r-1) g_{r-1} &\\
\ \     =(r-1)(g_{r-2}+\lambda g_{r-1}), & {\rm\ when\ } 
    3\leq r\leq n+1, \\
(n+1)g_n, &{\rm\ when\ } r=n+2.
\end{array} \right . 
\label{eq:com6}
\end{equation}
Recall that $\sha_C(X)$ is a free $C[x]$-module on the 
set $\calx$. 
So combining equation~(\ref{eq:com5}) and (\ref{eq:com6}), 
we obtain 

\begin{eqnarray*}
\lambda g_1  =  1, & \\
(r-1)(g_{r-2}+\lambda g_{r-1})=0, & 3\leq r\leq n+1, \\
(n+1)g_n = 0. 
\end{eqnarray*}
Since the characteristic of $C$ is zero by assumption, 
we have 
\begin{eqnarray*}
\lambda g_1  =  1, & \\
g_{r-2}+\lambda g_{r-1}=0, & 3\leq r\leq n+1, \\
g_n = 0. 
\end{eqnarray*}
Thus we have 
$g_n=0$, $g_{n-1}=-\lambda g_n=0, \ldots, 
g_1=-\lambda g_2 =0$. 
This contradicts with $\lambda g_1=1$,  
proving Theorem~\ref{thm:nacc1} 
when $X=\{x\}$. 

Let $X$ be any non-empty set. Fix an element $x_0\in X$. 
Then the surjective map   
$X\to \{x_0\}$ sending all $x\in X$ to $x_0$ 
induces a surjective homomorphism 
$\sha_C(X) \to \sha_C(\{x_0\})$ of Baxter algebras. 
In fact, the homomorphism $\sha_C(\{x_0\})\to \sha_C(X)$ 
induced by $\{x_0\}\to X,\ x_0\mapsto x_0$ provides a 
section of the first homomorphism. 
If $\sha_C(X)$ were a noetherian, then its surjective image 
$\sha_C(\{x_0\})$ would have to be noetherian also. 
We have already shown above that this is impossible. 
So $\sha_C(X)$ is not noetherian. 
\proofend

\section{Ascending chain condition for Baxter ideals}
\margin{sec:bacc}
\label{sec:bacc}

We now consider $\sha_C(C)$ in the category of Baxter algebras.
We first give some definitions. 

\begin{defn}
\begin{enumerate}
\item
A Baxter algebra $(R,P)$ is called a {\bf noetherian
Baxter algebra} if the set of Baxter ideals of
$(R,P)$ satisfies the ascending chain condition.
\item
A Baxter ideal $I$ of $(R,P)$
is called {\bf Baxter finitely generated} if there are finitely 
many elements $f_1,\ldots,f_r$ of $R$ such that $I$ is
the smallest Baxter ideal of $R$
containing $f_1,\ldots,f_r$.
\end{enumerate}
\end{defn}

\subsection{The case when $A=C$}

\begin{theorem}
\margin{thm:acc2}
\label{thm:acc2}
If $C$ is a noetherian ring, then
$\sha_C(C)$ and $\widehat{\sha}_C(C)$ are
noetherian Baxter algebras.
\end{theorem}

\begin{coro}
If $C$ is a noetherian ring, then any irreducible
Baxter $C$-algebra is a noetherian $C$-algebra. 
\end{coro}
\proof
This follows from Theorem~\ref{thm:acc2} since any irreducible 
Baxter $C$-algebra is a quotient of the free Baxter 
algebra $\sha_C(C)$.
\proofend

\noindent
{\bf Proof of Theorem~\ref{thm:acc2}: }
It is easy to see that $R$ is a noetherian Baxter
algebra if and only if every Baxter ideal $I$ of $R$ is
Baxter finitely generated. 
So we only need to prove that any Baxter ideal of $\sha_C(C)$
is Baxter finitely generated.
The idea of the proof is the same as that of
Theorem~\ref{thm:acc1}, following the Hilbert basis theorem,
except that multiplying by $x$ is
replaced by applying the Baxter operator $P_C$.
Let $I\subseteq \sha_C(C)$ be an Baxter ideal.
For each $j\in\NN$, let
\[\Sigma_j =\{ b_j\in C\mid \exists f_j\in I,\ {\rm such\ that\ }
    f_j=\sum_{k=0}^j b_k \bfone_k\}. \]
Then the same argument as in Theorem~\ref{thm:acc1}
shows that 
$\Sigma_j$ is an ideal of $C$.
Also
$b_j\in \Sigma_j$ implies that there 
exists a $f_j\in I$ such that 
$f_j=\sum_{k=0}^j b_k \bfone_k$. 
Then $P_C(f_j)=\sum_{k=0}^j b_k \bfone_{k+1}$
and so $b_j$ is in $\Sigma_{j+1}$. 
Thus $\Sigma_j\subseteq \Sigma_{j+1}$.
Since $C$ is noetherian, this chain of ideals stabilizes,
say at $j=m$.
Then $\Sigma_m=\cup_{j=1}^\infty \Sigma_j$,
and is finitely generated. 
Then as in the proof of Theorem~\ref{thm:acc1}, we construct from
this a set of elements of $I$ 

\[\{f_1^{(0)},\ldots,f_{k_0}^{(0)},f_1^{(1)},\ldots,f_{k_1}^{(1)},
\ldots, f_{1}^{(m)},\ldots,f_{k_m}^{(m)}\}\]
and prove that $I$ is the Baxter ideal generated by this set.

The statement for $\widehat{\sha}_C(C)$ can be proved in
the same way, replacing $\Sigma_j$ by
\[\Omega_j =\{ b_j\in C\mid \exists f_j\in I,\ {\rm such\ that\ }
    f_j=\sum_{k=j}^\infty b_k \bfone_k\}. 
\ \ \     \mbox{\proofend} \]

\subsection{The general case}
Because of Theorem~\ref{thm:acc2}, it is natural to ask
whether other Baxter algebras are noetherian, 
and in particular, whether the free Baxter algebras are
noetherian. Theorem~\ref{thm:acc2} shows that
if $C$ is noetherian, then $\sha_C(X)$ is noetherian if
$X$ is the empty set.
We will prove a theorem on free Baxter algebra 
$\sha_C(A)$. Consequences on $\sha_C(X)$ will be given 
in the corollary. 

Consider the $A$-module $A\otimes A$ with $A$ acting on the left 
tensor factor.
Let $M$ be a $C$-submodule of $A$, we use 
$\overline{A\otimes M}$ to denote the subgroup of $A\otimes A$ 
generated by elements of the form $a\otimes b,\ a\in A,\ b\in M$. 
It is a $A$-submodule of $A\otimes A$. 
It is easy to see that $\overline{A\otimes M}$ is 
the image of $A\otimes M$ in $A\otimes A$ under the 
natural map $A\otimes M\to A\otimes A$ induced by 
$M\hookrightarrow A$.

\begin{theorem}
\margin{thm:nacc2}
\label{thm:nacc2}
Let $A$ be a $C$-algebra. 
\begin{enumerate}
\item
Let $\calm$ be the partially ordered set consisting of 
$A$-submodules of $A\otimes A$ of the form 
$\overline{A\otimes M}$ where $M$ runs through $C$-submodule 
of $A$. 
If $\calm$ does not satisfy the ascending chain condition, 
then $\sha_C(A)$ of weight zero does not satisfy 
the ascending chain condition for Baxter ideals. 
\item
If $A$ is not a noetherian ring, 
then $\sha_C(A)$ of any weight does not satisfy 
the ascending chain condition for Baxter ideals. 
\end{enumerate}
\end{theorem}

\begin{remark}
The condition in (1) implies that $A$ is not a 
noetherian $C$-module. 
The following example shows that if the condition in (1) is 
weakened to the condition that $A$ is not a noetherian $C$-module, 
then the conclusion of (1) might not hold. 
To see what extra restriction is needed, see (1) in 
Corollary~\ref{co:nacc2}. 
\end{remark}
\begin{example} 
Let $C=\ZZ$ and $A=\QQ$. An infinite ascending chain of 
$\ZZ$-modules of $\QQ$ is given by 
\[ \ZZ \subsetneq 2^{-1}\ZZ \subsetneq \ldots \subsetneq 2^{-n}\ZZ 
    \subsetneq \ldots. \]
So the $\ZZ$-submodules of $\QQ$ does not satisfy the ascending chain condition. 
On the other hand, it is easy to verify that 
$\QQ\otimes \QQ=\QQ\otimes 1$ and it is the only 
$\QQ$-submodule of $\QQ\otimes \QQ$. 
Thus $\calm=\{\QQ\otimes \QQ\}$, trivially satisfying the ascending chain condition. 
Now for each $n\geq 1$, 
\[\QQ^{\otimes n} =\QQ \bfone^{\otimes n}\cong 
\underbrace{\QQ\otimes_\QQ \ldots \otimes_\QQ \QQ}_{n-{\rm factors}}\] 
as $\QQ$-modules. This implies that $\sha_{\ZZ}(\QQ)$ is isomorphic 
to $\sha_{\QQ}(\QQ)$ as rings. Since 
$\sha_{\QQ}(\QQ)$ is a noetherian ring by Theorem~\ref{thm:acc1}, 
$\sha_{\ZZ}(\QQ)$ is a noetherian ring. In particular, it has the ascending chain condition 
for Baxter ideals. 
\end{example}

Applying Theorem~\ref{thm:nacc2} to the case when $A=C[X]$, we obtain 
\begin{coro}
\margin{co:nacc2}
\label{co:nacc2}
\begin{enumerate}
\item
Let $A$ be a $C$-algebra. If $A$ is not a noetherian $C$-module, 
and if there is a $C$-linear homomorphism $A\to C$, 
then $\sha_C(X)$ of weight zero does not have the ascending chain condition for 
Baxter ideals. 
\item
If $X$ is not empty, then $\sha_C(X)$ of 
weight zero does not have the ascending chain condition for Baxter ideals.
\item
If $X$ is infinite, then $\sha_C(X)$ of any weight 
does not have the ascending chain condition for Baxter ideals. 
\end{enumerate}
\end{coro}

\proof
(1)
Denote the $C$-linear homomorphism $A\to C$ by $f$. 
By its $C$-linearity, $f$ must be surjective. 
Since the tensor product functor is right exact, for any 
$C$-module $M$, the surjective $C$-linear map 
$f: A\to C$ induces surjective abelian group homomorphism 
$f\otimes \id_M: A\otimes M \to C\otimes M\cong M$. 
Since $A$ is not a noetherian $C$-module, there are $C$-modules 
$M_n,\ n\geq 1$ such that $M_n\subsetneq M_{n+1}$ for all $n$. 
Suppose $\overline{A\otimes M_n} =\overline{A\otimes M_{n+1}}$ 
for some $n$.  
Let $j_{n,n+1}:M_n\rightarrow M_{n+1}$ and $j_{n+1}: M_{n+1} 
\rightarrow A$ be the natural embeddings. 
Consider the commutative diagram 
\[\begin{CD}
A\otimes M_n @>>> C\otimes M_n @>>> M_n \\
@V \id_A\otimes j_{n,n+1}VV @V \id_C\otimes j_{n,n+1} VV 
    @VV j_{n,n+1}V \\
A\otimes M_{n+1} @>>> C\otimes M_{n+1} @>>> 
    M_{n+1}\\
@V \id_A\otimes j_{n+1}VV  @V \id_C\otimes j_{n+1}VV    
    @VV j_{n+1} V\\
A\otimes A @>>>  C\otimes A @>>> A
\end{CD}
\]
where all horizontal maps on the left column are surjective 
and all horizontal maps on the right column are isomorphisms.
{}From $\overline{A\otimes M_n}=\overline{A\otimes M_{n+1}}$ we have 
\begin{eqnarray*}
\lefteqn{ (\id_A\otimes j_{n+1})(A\otimes M_{n+1}) = 
\overline{A\otimes M_{n+1}} }\\
& =& \overline{A\otimes M_n} \\
&=& (\id_A\otimes j_n)(A\otimes M_n)\\
&=& 
((\id_A\otimes j_{n+1})\circ (\id_A\otimes j_{n,n+1})) (A\otimes M_n).
\end{eqnarray*}
Since all the horizontal maps in the commutative diagrams are surjective,  
we further have 
\[ M_{n+1}= j_{n+1}(M_{n+1}) = 
(j_{n+1}\circ j_{n,n+1}) (M_n) =M_n.
\]
This is a contradiction. 
Therefore the assumption in the first statement of 
Theorem~\ref{thm:nacc2} is satisfied, proving that 
$\sha_C(A)$ of weight zero does not have the ascending chain condition for Baxter ideals. 

\mvp
\noindent
(2)
Let $X$ be non-empty. We only need to show that $C[X]$ satisfies 
the assumptions in the first statement of the corollary, 
which follows from the first statement of the theorem. 

For each $n\geq 1$, let $M_n$ be the $C$-submodule of $C[X]$ 
generated by $x^k,\ 1\leq k\leq n$. 
Since, for each $n\geq 1$, $M_n$ is a submodule of $M_{n+1}$ and 
$x^{n+1}$ of $M_{n+1}$ is not in $M_n$, we have an infinite 
ascending chain $M_n\subsetneq M_{n+1}$. 
So $C[X]$ is not a noetherian $C$-module. 
On the other hand, 
the $C$-algebra map $f:C[X]\to C$ induced by sending $x\in X$ to 
$0$ is clearly the $C$-linear map we want. Therefore (1) applies. 

\mvp
\noindent
(3) This follows from the second statement of Theorem~\ref{thm:nacc2} 
since $C[X]$ is not a noetherian ring when $X$ is infinite. 
\proofend

Before the proof of Theorem~\ref{thm:nacc2}, we first prove a 
lemma. 

\begin{lemma}
\margin{lem:ideal}
\label{lem:ideal}
Let $M$ be a $C$-submodule of $A$. Define
\[ S =\bigcup_{k=1}^\infty
    \left \{ \otimes_{i=0}^k x_i\in A^{\otimes (k+1)} \mid
x_{i_0}\in M {\rm \ for\ some\ } 1\leq i_0\leq k
    \right \}.\]
Let $I$ be the abelian subgroup of $\sha_C(A)$ generated
by $S$.
If either the weight of $\sha_C(A)$ is zero or $M$ is an ideal of $A$, 
then $I$ is a Baxter ideal of $\sha_C(X)$. 
In fact, $I=I\tp$, the Baxter ideal of $\sha_C(A)$ generated by 
$P_A(M)$. 
\end{lemma}
\proof
We only need to prove the last statement. 
We first prove that $I\subseteq I\tp$. 
For this we only need to show that, for each 
$x=\otimes_{i=0}^k x_i\in S$, we have $x \in I\tp$.
We prove by induction on $k\geq 1$.
When $k=1$, we have $x=x_0\otimes x_1$ with $x_1\in M$.
Then $x=x_0 (1\otimes x_1)\in I_n\tp$.
Now let $x=\otimes_{i=0}^{k+1}x_i\in S$, so $x_i\in A$
and $x_{i_0}\in M$ for some $1\leq i_0\leq k+1$.
If $i_0>1$, then
$x_1\otimes \ldots x_{k+1}$ is in $S$,
and hence by the induction hypothesis, is in $I\tp$.
Then $x=x_0 P_A(x_1\otimes\ldots\otimes x_{k+1})$
is in $I\tp$. 
If $i_0=1$,
consider the equation obtained by the definition of the 
mixable shuffle product~(\ref{eq:shuf})
\begin{eqnarray*}
\lefteqn{
 (x_0\otimes x_2\otimes \ldots \otimes x_{k+1})
     (\bfone_A\otimes x_1)
=x_0\otimes x_1\otimes\ldots\otimes x_{k+1}}\\
&+& \sum_{j=2}^{k} x_0\otimes x_2\otimes
    \ldots\otimes x_j\otimes x_1 \otimes
    x_{j+1}\ldots\otimes x_{k+1} \\
&+& x_0\otimes x_2\otimes \ldots \otimes x_{k+1}\otimes x_1\\
&+& \lambda \sum_{j=2}^{k+1} x_0\otimes x_2\otimes
    \ldots\otimes x_j x_1 \otimes
    x_{j+1}\ldots\otimes x_{k+1}.
\end{eqnarray*}
Since $x_1$ is in $M$, we see that  
$\bfone_A\otimes x_1$ is in $I\tp$. So the left hand side of the
equation is in $I\tp$.
Again because $x_1$ is in $M$, 
the induction hypothesis implies that 
every term on the right hand side except the first term
and the terms in the last sum are in $I\tp$. 
If $\lambda=0$, then the last sum disappears.  
So the first term is also in $I\tp$.
On the other hand, if $M$ is an ideal of $A$, then 
$x_j x_1$ is in $M$ for $2\leq j\leq k+1$. 
Hence by induction hypothesis, every term in the last sum is in 
$I\tp$. So again the first term is in $I\tp$. 
This completes the induction, proving that 
$I\subseteq I\tp$. 

We next prove that $I$ contains $I\tp$. 
For this we only need to show that $I$ is a Baxter ideal of $\sha_C(A)$ 
since $I$ clearly contains $M$. 
By the definition of $S$ we have $P_A(S)\subseteq S$. 
So we get $P_A(I)\subseteq I$.
Since $I$ is clearly closed under addition, 
it remains to verify that,
if $x\in I$ and $y\in \sha_C(A)$, then 
$x y$ is in $I$. For this we only need to verify this property
for $x=\otimes_{i=0}^k x_i\in S $ and
$y=\otimes_{j=0}^m y_j \in A^{\otimes (m+1)},\ m\geq 0$. 
By definition, 
\margin{eq:shuf2}
\begin{equation}
x\ y
= x_0y_0\otimes \sum_{(\sigma,T)\in \bs(m,n)} \lambda^{\mid T\mid}
    \sigma(x_1\otimes \ldots \otimes x_m\otimes y_1
    \otimes \ldots \otimes y_n;T).
\label{eq:shuf2}
\end{equation}
For each $(\sigma,T)\in \bs(m,n)$, the set of $(m,n)$-mixable 
shuffles defined in \S~\ref{sec:shuf}, 
write
\[ \sigma(x_1\otimes \ldots \otimes x_m\otimes y_1
    \otimes \ldots \otimes y_n;T)
    =z_1\hts\ldots \hts z_{m+n}. \]
Then
$(z_1,\ldots,z_{m+n})$ is a permutation of
$(x_1,\ldots,x_m,y_1,\ldots,y_n)$. 
Since $x_{i_0}$ is in $M$ for some $1\leq i_0\leq n$,
one of $z_i,\ 1\leq i\leq m+n$, is in $M$.
If $\lambda=0$, then the only non-zero terms in the sum on the right 
hand side of equation~(\ref{eq:shuf2}) are of the form 
$z_1\otimes \ldots \otimes z_{m+n}$. 
So the right hand side is in $I$. 
Now assume that $M$ is an ideal of $A$.  
For any $(\sigma,T)\in \bar{S}(m,n)$, by the definition of 
$\sigma(x_1\otimes\ldots\otimes x_m\otimes y_1\otimes\ldots\otimes y_n;T)$, 
either $x_{i_0}$ or $x_{i_0}y_j$ for some $1\leq j\leq n$ 
occurs as one of the tensor factors in 
\[\sigma(x_1\otimes\ldots\otimes x_m\otimes y_1\otimes\ldots\otimes y_n;T)
=z_1\hts\ldots \hts z_{m+n}. 
\] 
Since $x_{i_0}$ and $x_{i_0} y_j$ are in $M$, we see that 
\[\lambda^{\mid T\mid} x_0y_0\otimes 
\sigma(x_1\otimes\ldots\otimes x_m\otimes y_1\otimes\ldots\otimes 
y_n;T)\]
is in $S$. 
Thus $x y$ is in $I$.
Therefore $I$ is an Baxter ideal of $\sha_C(A)$. 
Consequently, $I\tp\subseteq I$. 
The lemma is now proved.
\proofend

\noindent
{\bf Proof of Theorem~\ref{thm:nacc2}:}
(1) By assumption, 
there are $C$-submodules $M_n,\ n\in \NN_+$,  of $A$
such that $\overline{A\otimes M_n}$ is a proper submodule of 
$\overline{A\otimes M_{n+1}}$ for all $n$. 
Define
\[ S_n=\bigcup_{k=1}^\infty
    \left \{ \otimes_{i=0}^k x_i\in A^{\otimes (k+1)} \mid
x_{i_0}\in M_n {\rm \ for\ some\ } 1\leq i_0\leq k
    \right \}.\]
Let $I_n$ be the abelian subgroup of $\sha_C(A)$ generated
by $S_n$.
Since we assume that $\lambda$ is zero, 
by Lemma~\ref{lem:ideal}, $I_n$ is a Baxter ideal of $\sha_C(X)$. 
Suppose $\sha_C(X)$ satisfies the ascending chain condition for Baxter ideals. 
Then the ascending chain $I_n,\ n\geq 1$ stabilizes for large $n$. 
In particular, $I_n=I_{n+1}$ for some $n$. 
Let 
\[ p:\sha_C(A)=\oplus_{m=1}^\infty A^{\otimes m} 
\rar A\otimes A\]  
be the natural projection from $\sha_C(A)$ onto the summand with $m=2$.
Then from the construction of $I_n$ we have
$p(I_n)=\overline{A\otimes M_n}$.
Thus $I_n=I_{n+1}$ implies that 
$\overline{A\otimes M_n}=\overline{A\otimes M_{n+1}}$. 
This contradicts the choice of $M_n$. 
So $\sha_C(A)$ does not have the ascending chain condition for Baxter ideals. 

\mvp

The proof of (2) is similar. 
\proofend

\begin{center}
{\bf ACKNOWLEDGEMENTS }
\end{center}
The author thanks William Keigher for helpful discussions.

\addcontentsline{toc}{section}{\numberline {}References}

\end{document}